# Evaluation for moments of a ratio with application to regression estimation

PAUL DOUKHAN[1] and GABRIEL LANG[2]

[1]*UMR 8088 Analyse, Géométrie et modélisation, Université Cergy-Pontoise, 2, avenue Adolphe Chauvin, 95302 Cergy-Pontoise Cedex, France.*
*E-mail: paul.doukhan@free.fr*
[2]*UMR 518 Mathématique et Informatique appliquées, AgroParisTech ENGREF, 19 avenue du Maine, 75732 Paris Cedex 15, France.*
*E-mail: gabriel.lang@engref.agroparistech.fr*

Ratios of random variables often appear in probability and statistical applications. We aim to approximate the moments of such ratios under several dependence assumptions. Extending the ideas in Collomb [*C. R. Acad. Sci. Paris* **285** (1977) 289–292], we propose sharper bounds for the moments of randomly weighted sums and for the $L^p$-deviations from the asymptotic normal law when the central limit theorem holds. We indicate suitable applications in finance and censored data analysis and focus on the applications in the field of functional estimation.

## 1. Introduction

We consider statistics with the form of a ratio. This situation arises naturally in the following simple and generic example. Let $(V_i, W_i)_{i \geq 0}$ be a stationary sequence with values in a finite space $\mathcal{V} \times \mathcal{W}$, with $\mathcal{V} \subset \mathbb{R}$. A conditional expectation can be expressed as

$$\mathbb{E}(V_0|W_0 = w) = \frac{\sum_{v \in \mathcal{V}} v \mathbb{P}(\{v, w\})}{\sum_{v \in \mathcal{V}} \mathbb{P}(\{v, w\})}.$$

Two examples of statistical ratios can be derived from this expression:

- For discrete random variables, it is empirically estimated from a sample $(V_i, W_i)_{1 \leq i \leq n}$ by the random quantity

$$\widehat{\mathbb{E}}(V_0|W_0 = w) = \frac{1/n \sum_{i=1}^n V_i \cdot \mathbb{1}\{W_i = w\}}{1/n \sum_{i=1}^n \mathbb{1}\{W_i = w\}}.$$

- The case of real-valued data is more involved as $\mathbb{P}(\{v,w\})$ has no rigorous meaning, but standard smoothing techniques allow us to consider extensions. Replacing $\mathbb{1}\{W_i = w\}$ by an approximate Dirac measure $\delta_n(w, W_i)$, the estimate corresponds to the Nadaraya–Watson kernel estimator that will be studied in detail in the sequel.

In this paper, we consider a ratio of two empirical quantities, namely

$$\widehat{R}_n = \frac{\widehat{N}_n}{\widehat{D}_n}, \qquad \widehat{D}_n = \frac{1}{n}\sum_{i=1}^n U_{i,n}, \widehat{N}_n = \frac{1}{n}\sum_{i=1}^n U_{i,n}V_{i,n}, \tag{1}$$

$U_{i,n}$ and $V_{i,n}$ being two arrays of random variables. Examples of this are:

- Functional estimation of a conditional expectation: let $(X_i, Y_i) \in \mathbb{R}^d \times \mathbb{R}$ be a stationary process and $K$ a kernel function. If we define

$$U_{i,n} = K((X_i - x)/h_n)/h_n^d, \qquad V_{i,n} = Y_i,$$

then $\widehat{R}_n = \widehat{r}(x)$ is an estimator for $r(x) = \mathbb{E}(Y_i|X_i = x)$ and $h_n \to 0$, $nh_n^d \to \infty$ as $n \to \infty$; see Tsybakov [28] for a general setting and Ango Nze and Doukhan [3] for dependent data cases.
- Computation of empirical means for censored data: let the censoring $U_i = C_i \in \{0, 1\}$ be independent of a process $(V_i)$ and assume that $V_i$ is observed if and only if $C_i = 1$.

$$\widehat{R}_n = \frac{1}{\#\{i \in \{1,n\}|C_i = 1\}} \sum_{1 \leq i \leq n, C_i = 1} V_i.$$

A example of this situation is the estimation of covariances of a process $X$ under censoring where $V_i = X_i X_{i+\ell}$. Under stationarity, the covariance function is $\gamma_X(\ell) = \gamma_Y(\ell)/(\gamma_C(\ell) + \mathbb{E}C_0^2)$, where $Y_i = C_i X_i$ is observed. Furthermore, moments of the empirical covariances are used to build the periodogram from the censored data.
- General weighted sums

$$\widehat{R}_n = \left(\sum_{1 \leq i \leq n} U_i V_i\right) \Big/ \left(\sum_{1 \leq i \leq n} U_i\right)$$

may be used to model various quantities like prices, with prices per unit $V_i$ and volumes $U_i$, as in [20].

*Various alternative questions also involve a division:*

- Functional estimation of point processes. A compound Poisson processes (CPP) can be expressed as $\xi = \sum_{j=1}^N \alpha_j \delta_{X_j}$ for some Poisson variable $N$ and some random process $(\alpha_j, X_j)_{j \geq 1}$, $\alpha_j > 0$, $X_j \in R^d$. For a sequence of mixing couples of CPP $(\xi_i, \eta_i)_{i \geq 1}$ with $\mu = \mathbb{E}\eta_1 \ll \nu = \mathbb{E}\xi_1$, Bensaïd and Fabre [5] estimate the Radon–Nikodym density $\varphi = d\mu/d\nu$ with kernel estimates $\varphi_n = g_n/f_n$, with $f_n(x) = \sum_{i=1}^n \eta_i \star K_n(x)$, $g_n(x) = \sum_{i=1}^n \xi_i \star K_n(x)$, where, for example, $\xi \star K_n(x) = \sum_{j=1}^N \alpha_j K_n(X_j - x)$ with $K_n(x) = K(x/h)/h^d$. The procedure is thus analogous to the Nadaraya–Watson estimator. Quadratic errors of this ratio are bounded under the assumption $|\varphi_n(x)| \leq C$ that can be easily relaxed using our result.
- Self-normalized sums, for example, in [11].
- Simulation of Markov chains and the Monte Carlo MC technique widely developed in the monograph of Robert and Casella [23]; a more precise reference is Li and Rabitz [21]; see relations (45)–(47), which explicitly involve ratios for reducing the dimensionality in a nonparametric problem.
- Particle filtering, considered from the theoretical viewpoint in [12], and for applications to change point problems in [18].

Deducing the convergence in probability of the ratio from the convergence of the denominator and numerator is straightforward, but in some statistical problems, $\mathbb{L}^p$-convergence has to be checked. Evaluating the moments of $\widehat{R}_n$ is much more difficult, even if one knows sharp bounds of moments for both the numerator and the denominator. Curiously, we did not find many references on this subject. One method is to compute the exact distribution of the ratio, as Spiegelmann and Sachs [25] did for the moments of a Nadaraya–Watson regression estimator with $\{0,1\}$-valued kernels. In this case, independence allows the use of binomial-based distributions. However, such computations are generally difficult to handle. An alternative is the expansion in [6]. We addressed this problem for a dependent data frame in the paper [7], published after Gérard Collomb's death. In [6] and [7], Collomb assumed that convergence rates in $\mathbb{L}^q$ for $q > 2p$ are known for the denominator.

This limitation is avoided here by using an interpolation technique and we shall only assume such rates for some $q > p$. With the notation (1), we set

$$N_n = \mathbb{E}\widehat{N}_n, \qquad D_n = \mathbb{E}\widehat{D}_n \quad \text{and} \quad R_n = \frac{N_n}{D_n}. \tag{2}$$

We aim to provide $\mathbb{L}^p$-rates of convergence to 0 of the expression

$$\Delta_n = \widehat{R}_n - R_n. \tag{3}$$

In some of our applications, the expectations $N_n$ and $D_n$ are constant. In other cases, they converge to some constants $N$ and $D$ as $n \to \infty$, and the moments of the ratio may be proven to converge with the bound

$$\left\| \widehat{R}_n - \frac{N}{D} \right\|_p \leq \|\Delta_n\|_p + \left| \frac{N}{D} - R_n \right|.$$

Convergence in probability or a.s. is immediate, but to obtain moment bounds, one has to divide by a non-zero expression; for simplicity, from now on, we will assume that $U_{i,n} \geq 0$. The previous expression is then also a weighted sum

$$\widehat{R}_n = \sum_{i=1}^n w_{i,n} V_{i,N}, \qquad w_{i,n} = \frac{U_{i,n}}{\sum_{j=1}^n U_{j,n}} \geq 0, \qquad \sum_{i=1}^n w_{i,n} = 1$$

so that $\widehat{R}_n$ belongs to the convex hull of $(V_{i,n})_{1 \leq i \leq n}$.

The paper is organized as follows. Section 2 is devoted to the main lemma and comments. The two following sections are dedicated to its applications: to simple weighted sums in Section 3 and to the Nadaraya–Watson kernel regression estimation $\widehat{r}_n(x)$ of a regression function $r(x) = \mathbb{E}(Y|X=x)$ in Section 4. The latter is divided into two subsections: the first subsection directly applies the lemma to provide the minimax bound $\|\widehat{r}_n(x) - r(x)\|_p = \mathcal{O}(n^{-\rho/(2\rho+d)})$ in an estimation problem with dimension $d$ and regularity $\rho$ less than 3; the second subsection makes use of the same ideas with a slight modification to derive the (also minimax) bound $\|\sup_{x \in B} |\widehat{r}_n(x) - r(x)|\|_p = \mathcal{O}((n/\log n)^{-\rho/(2\rho+d)})$ over a suitable compact subset $B \subset \mathbb{R}^d$ (see Stone [26] and [27]). We prove our result under various dependence assumptions: independent, strongly mixing, absolutely regular or weakly dependent (either causal or non-causal) sequences. The last section includes the proofs; it consists of four subsections, devoted to the main lemma, weighted sums, moments of Nadaraya–Watson estimation and sup bounds of this estimator, respectively.

## 2. Main lemmas

Lemma 1 means that for $q$ slightly larger than $p$, the rate of the $q$th order moment of the denominator and of the $p$th order moment of the numerator allow us to derive a bound of the rate for the $p$th order moment of the ratio.

**Lemma 1.** *Assume that $\|\widehat{N}_n - N_n\|_p \leq v_n$ and $\|\widehat{D}_n - D_n\|_q \leq v_n$ for some $q > p$. If, moreover, $\|U_{i,n} V_{i,n}\|_r \leq C_n$ and $\|V_{i,n}\|_s \leq c_n$, where $q/p - q/r \geq 1$, $1/p > 1/q + 1/s$, then*

$$D_n \|\Delta_n\|_p \leq \left(1 + \frac{|N_n|}{D_n} + \frac{|N_n|^\beta v_n^{1-\beta}}{D_n} + \frac{C_n^\beta v_n^{1-\beta}}{D_n} + \frac{v_n^\alpha c_n n^{1/s}}{D_n^\alpha}\right) v_n,$$

*where $\alpha, \beta$ are chosen from the parameters $p, q, r, s$ by setting*

$$\alpha = q\left(\frac{1}{p} - \frac{1}{s} - \frac{1}{q}\right) \leq 1 \leq \beta^{-1} = q\left(\frac{1}{p} - \frac{1}{r}\right).$$

*Remarks.*

- In all of our examples, $c_n \equiv c$ will be a constant.
- If $C_n \equiv C$ is also a constant, then we assume that $r = \frac{pq}{q-p}$ so that $\beta = 1$. In this case, large values of $q$ give $r$ close to (and larger than) $p$; if now $q > p$ is very close to $p$, then $r$ needs to be very large and $s$ even larger.

  *This is the situation for weighted sums or censored data questions.*

  Here, $V_{i,n} = V_i$ and $0 \leq U_{i,n} = U_i$, and the sequence $(U_i, V_i)$ is stationary. Moreover, $v_n = c/\sqrt{n}$ and thus $\Delta_n = \mathcal{O}(n^{-1/2})$ if $\alpha = 2/s$. This condition can be expressed as $s = p(q+2)/(q-p)$, as proved in the forthcoming Theorem 1.

- When the sequence $C_n$ is not bounded, in order to control the corresponding term, we shall use an exponent $\beta < 1$ and assume that $1/p > 1/q + 1/r$. The order $q$ of the moment of the denominator should be larger, as well as the order of the moments of the variables $V_{i,n}$ (in the case of functional estimation, for example). Here, $v_n = c/\sqrt{nh_n^d} \gg 1/\sqrt{n}$ as $h_n \to_{n\to\infty} 0$ and $nh_n^d \to_{n\to\infty} \infty$.
- *Orlicz spaces.* Instead of $\mathbb{L}^q$-norms, we may consider Orlicz norms and ask only for $x^p \log^q x$-order moments of the denominator. Exponential moments of the variables $V_{i,n}$ would be used because of the relations (4), (14) and of the Pisier inequalities (15).
- *Suprema.* The same equations (4), (14) and (15) are adapted to derive bounds of suprema for moments for expressions involving an additional parameter; an emblematic example of this situation is the regression estimation given in Section 4.2.

We consider two distinct classes of applications in Sections 3 and 4, devoted respectively to *weighted sums* and *nonparametric regression.* The following inequality is essential to bound the uniform rates of convergence of a Nadaraya–Watson regression estimator and it is thus presented as a specific lemma.

**Lemma 2.** *Letting* $0 < \alpha < 1$, *we have*

$$D_n|\Delta_n| \leq |\widehat{N}_n - N_n| + \frac{|\widehat{N}_n|}{D_n}|\widehat{D}_n - D_n| + \max_{1 \leq i \leq n} |V_{i,n}| \frac{|\widehat{D}_n - D_n|^{1+\alpha}}{|D_n|^\alpha}. \tag{4}$$

Inequality (4) also implies tail bounds for $\Delta_n$'s distribution.

$\mathbb{L}^{p'}$-convergence may also be addressed as follows, as suggested by an anonymous referee.

**Corollary 1.** *Assume that* $v_n^{-1}\Delta_n \to_{n\to\infty} Z$ *converges in distribution to some* $Z$ *such that* $\|Z\|_{p'} < \infty$ *for some* $p' < p$. *If the conditions in Lemma 1 hold and the bound is such that* $\|\Delta_n\|_p \leq cv_n$ *for some* $c > 0$, *then*

$$v_n^{-1}\|\Delta_n\|_{p'} \underset{n\to\infty}{\longrightarrow} \|Z\|_{p'}.$$

*Remarks.*

- For weighted sums, a central limit theorem is obtained in all the cases considered below so that the result applies.
- For the Nadaraya–Watson estimator, Ango Nze and Doukhan [1] prove that $\sqrt{nh_n^d}(\widehat{r}(x) - r(x)) \to N(a(x), b(x))$ for functions $a(x), b(x)$. For the case of bounded regressors, the result also holds for weakly dependent cases (see Doukhan and Louhichi [14]).

## 3. Weighted sums

We consider here the simplest application of Lemma 1. Let $(U_i, V_i)_{i \in \mathbb{Z}}$ be a stationary sequence and set $\widehat{D}_n = \sum_{i=1}^n U_i/n$, $\widehat{N}_n = \sum_{i=1}^n U_i V_i/n$. Then $N_n = N = \mathbb{E}U_1 V_1$, $D_n = D = \mathbb{E}U_1$ and $\widehat{R}_n = \widehat{N}_n/\widehat{D}_n$, $R_n = R = N/D$.

**Theorem 1.** *Let $(U_i, V_i)_{i \in \mathbb{Z}}$ be a stationary sequence with $U_i \geq 0$ (a.s.). Let $0 < p < q$ and assume that $\|U_i V_i\|_r \leq c$ for $r = pq/(q-p)$ and $\|V_i\|_s \leq c$ for $s = p(q+2)/(q-p)$. If the dependence of the sequence $(U_i, V_i)_{i \in \mathbb{Z}}$ is such that*

$$\|\widehat{D}_n - D\|_q \leq Cn^{-1/2}, \qquad \|\widehat{N}_n - N\|_p \leq Cn^{-1/2}, \tag{5}$$

*then $\|\widehat{R}_n - R\|_p = \mathcal{O}(n^{-1/2})$.*

From now on, we assume that $\|V_i\|_s \leq c$ and $\|U_i V_i\|_r \leq c$, and prove that (5) holds.

### 3.1. Independent case

Assume that $(U_i, V_i)$ is i.i.d. Assume that $\|U_0\|_q \leq c$ and $\|U_0 V_0\|_p \leq c$. From the Marcinkiewikz–Zygmund inequality for independent variables,

$$\mathbb{E}|\widehat{D}_n - D|^q \leq C_q \mathbb{E}|U_1|^q n^{-q/2} \leq Cn^{-q/2},$$
$$\mathbb{E}|\widehat{N}_n - N|^p \leq C_p \mathbb{E}|U_1 V_1|^p n^{-p/2} \leq Cn^{-p/2}.$$

The Hölder inequality implies that the assumptions hold if $\|U_0\|_q$, and $\|V_0\|_{qp/(q-p)}$ are bounded.

### 3.2. Strong mixing case

Denote by $(\alpha_i)_{i \in \mathbb{N}}$ the strong mixing coefficients of the stationary sequence $(U_i, V_i)_{i \in \mathbb{N}}$.

**Proposition 1.** *Assume that for $s \geq r' > q$, $\|U_0\|_{r'} \leq c$. Relation (5) holds if $\alpha_i = \mathcal{O}(i^{-\alpha})$ with*

$$\alpha > \left(\frac{p}{2} \cdot \frac{r}{r-p}\right) \vee \left(\frac{q}{2} \cdot \frac{r'}{r'-q}\right).$$

### 3.3. Causal weak dependence

Let $(W_i)_{i \in \mathbb{N}}$ be a centered sequence with values in $\mathbb{R}^d$; for $k \geq 0$, we set $\mathcal{M}_k = \sigma(W_j;\ 0 \leq j \leq k)$. For each $i \in \mathbb{N}$, we define the $\gamma$ coefficients by

$$\gamma_i = \sup_{k \geq 0} \|\mathbb{E}(W_{i+k}|\mathcal{M}_k)\|_1.$$

**Proposition 2.** *Assume that $U_i$ and $V_i$ are two stationary $\gamma$-dependent sequences with common $\gamma_i = \mathcal{O}(i^{-\gamma})$. Assume that $U$ and $V$ are independent and that $\|U_0\|_\infty \leq c$. Relation (5) holds if*

$$\gamma > \left(\frac{p}{2} \cdot \frac{s-1}{s-p}\right) \vee \frac{q}{2}.$$

### 3.4. Non-causal weak dependence

Here, we consider non-causal weakly dependent stationary sequences and assume that $q$ and $p$ are integers. A sequence $(W_i)_{i \in \mathbb{N}}$ is $\lambda$-weakly dependent if there exists a sequence $(\lambda(i))_{i \in \mathbb{N}}$ decreasing to zero such that

$$\begin{aligned}|\operatorname{Cov}(g_1(W_{i_1},\ldots,W_{i_u}), g_2(W_{j_1},\ldots,W_{j_v}))| \\ \leq (u \operatorname{Lip} g_1 + v \operatorname{Lip} g_2 + uv \operatorname{Lip} g_1 \operatorname{Lip} g_2) \lambda(k)\end{aligned} \quad (6)$$

for any $u$-tuple $(i_1, \ldots, i_u)$ and any $v$-tuple $(j_1, \ldots, j_v)$ with $i_1 \leq \cdots \leq i_u < i_u + k \leq j_1 \leq \cdots \leq j_v$, where $g_1, g_2$ are real functions of $\Lambda^{(1)} = \{g_1 \in \Lambda | \|g_1\|_\infty \leq 1\}$ defined respectively on $\mathbb{R}^{Du}$ and $\mathbb{R}^{Dv}$ ($u, v \in \mathbb{N}^*$). Recall, here, that $\Lambda$ is the set of functions with $\operatorname{Lip} g_1 < \infty$ for some $u \geq 1$, with

$$\operatorname{Lip} g_1 = \sup_{(x_1,\ldots,x_u) \neq (y_1,\ldots,y_u)} \frac{|g_1(y_1,\ldots,y_u) - g_1(x_1,\ldots,x_u)|}{|y_1 - x_1| + \cdots + |y_u - x_u|}.$$

The monograph Dedecker *et al.* [9] details weak dependence concepts, models and results.

**Proposition 3.** *Assume that the stationary sequence $(U_i, V_i)_{i \in \mathbb{N}}$ is $\lambda$-weakly dependent for $\lambda_i = \mathcal{O}(i^{-\lambda})$. Assume that $p$ and $q \geq 2$ are even integers.*
*Relation (5) holds under each of the following sets of conditions:*

- *the processes $U$ and $V$ are independent, $\|U_0\|_\infty \leq c$ and $\lambda > \frac{q}{2}$;*
- *for $r' \leq s$, $\|U_0\|_{r'} \leq c$ and $\lambda > \frac{r'}{r'-2}\frac{q}{2}$.*

**Remark.**

- Non-integer moments $q \in (2,3)$ are considered in [17], Lemma 4, and the same inequality holds if $\mathbb{E}|Z_i|^{q'} < \infty$ with $q' = q + \delta$, and $\lambda(i) = \mathcal{O}(i^{-\lambda})$ with $\lambda > 4 + 2/q'$ for

$$q \leq 2 + \frac{1}{2}(\sqrt{(q'+4-2\lambda)^2 + 4(\lambda-4)(q'-2)-2} + q' + 4 - 2\lambda) \qquad (\leq q').$$

## 4. Regression estimation

We now use a measurable bounded function $K \colon \mathbb{R}^d \to \mathbb{R}$. Consider a stationary process $(X_i, Y_i) \in \mathbb{R}^d \times \mathbb{R}$. We now set, for some $h = h_n$,

$$U_{i,n} = U_{i,n}(x) = \frac{1}{h^d} K\left(\frac{X_i - x}{h}\right) \quad \text{and} \quad V_{i,n} = Y_i,$$

where $h$ tends to zero as $n$ tends to infinity. Then $\widehat{R}_n = \widehat{r}(x)$ is the Nadaraya–Watson estimator of $r(x) = \mathbb{E}(Y_i | X_i = x)$. Independently from the dependence structure of the process $(X_i, Y_i)$, we first introduce the following regularity conditions:

(A1) for the point of interest $x$, the functions $f, g$ are $k$-times continuously differentiable around $x$ and $(\rho - k)$-Hölderian, where $k < \rho$ is the largest possible integer;

(A2) the function $K$ is Lipschitz, admits a compact support, and satisfies $K(u) \geq 0$ ($\forall u \in \mathbb{R}^d$) and

$$\int_{\mathbb{R}^d} K(u)\,\mathrm{d}u = 1,$$

$$\int_{\mathbb{R}^d} u_1^{\ell_1} \cdots u_d^{\ell_d} K(u)\,\mathrm{d}u = 0 \qquad \text{if } 0 < \ell_1 + \cdots + \ell_d < k.$$

Moment and conditional moment conditions are also needed:

(A3) for the point $x$ of interest, there exist $r$ and $s$, with $r \leq s$ such that:
  1. $\|Y_0\|_s = c < \infty$;
  2. $g_r(x) = \int |y|^r f(x,y)\,\mathrm{d}y$ is a function bounded around the point $x$;
  3. $G(x,x) = \sup_i f_i(x,x)$ is bounded around the point $(x,x)$, where $f_i(x', x'')$ denotes the joint density of $(X_0, X_i)$.

*Remarks.*

- First, notice that the last condition holds immediately for independent sequences $(X_i)$ with a locally bounded marginal density.
- An alternative condition involving local uniform bounds of

$$(x', x'') \mapsto H(x', x'') = \sup_i \int |yy'| f_i(x', y'; x'', y'')\,\mathrm{d}y'\,\mathrm{d}y'',$$

where $f_i(x', y'; x'', y'')$ denotes the joint density of $(X_0, Y_0; X_i, Y_i)$, yields sharper results, but such conditions are generally difficult to check. They hold for independent sequences if $g_1$ is locally bounded.

We note here that the assumption $K \geq 0$ in (A2) implies that $k = 1$ or 2. We are not able to control biases by $h_n^\rho$ in Proposition 4 if $\rho \geq 3$. Without this non-negativity condition, the moments of numerator and denominator are still controlled, but we definitely cannot handle the moments of our ratio.

### 4.1. Moment estimation

We now consider the quantity $\delta_n(x) = \|\widehat{r}(x) - r(x)\|_p$ for $x \in \mathbb{R}^d$.

**Proposition 4 (Bias).** *Assuming that $f$ is bounded below around the point $x$ and that (A1) and (A2) hold, we have*

$$\left| r(x) - \frac{\mathbb{E}\widehat{g}(x)}{\mathbb{E}\widehat{f}(x)} \right| = \mathcal{O}(h_n^\rho) \qquad \forall x \in \mathbb{R}^d.$$

We now assume that the numerator and denominator satisfy the usual rate:

(A4) For the point of interest $x$ and for some $q > p$, for $h = h_n \to 0$ and $nh_n^d \to \infty$ as $n \to \infty$, there exists a constant $c > 0$ such that:
1. $\|\widehat{f}(x) - \mathbb{E}\widehat{f}(x)\|_q \le c/\sqrt{nh_n^d}$;
2. $\|\widehat{g}(x) - \mathbb{E}\widehat{g}(x)\|_p \le c/\sqrt{nh_n^d}$.

**Proposition 5.** *Assuming that $f$ is bounded below around the point $x$ and that assumptions (A2), (A3) and (A4) hold, we have*

$$\left\| \widehat{r}(x) - \frac{\mathbb{E}\widehat{g}(x)}{\mathbb{E}\widehat{f}(x)} \right\|_p \le C(1 + h^{d\beta(1/r-1)}(nh^d)^{(\beta-1)/2} + (nh^d)^{-\alpha/2} n^{1/s}) v_n$$

*with* $\alpha = q(\frac{1}{p} - \frac{1}{s} - \frac{1}{q})$ *and* $\beta = \frac{pr}{q(r-p)}$.

Those two propositions imply that the optimal window width $h \sim n^{-1/(2\rho+d)}$ equilibrates both expressions to get the minimax rate $\delta_n(x) = \mathcal{O}(n^{-\rho/(2\rho+d)})$:

**Theorem 2.** *Choose the window width $h_n = Cn^{-1/(2\rho+d)}$ for a constant $C > 0$. Assume that $f$ is bounded below around the point $x$ and that assumptions (A1), (A3) and (A4) hold for*

$$\frac{pd(r-1)}{qr - pq - pr} \vee \frac{pd}{qs - pq - ps - 2p} \le \rho. \qquad (7)$$

*There then exists a constant $C > 0$ such that $\|\widehat{r}(x) - r(x)\|_p \le Cn^{-\rho/(2\rho+d)}$.*

We now consider specific dependence structures to get the moment inequalities of (A4) for $\widehat{f}$ and $\widehat{g}$. From now on, fix $x$ and write

$$\widehat{g}(x) - \mathbb{E}\widehat{g}(x) = \frac{1}{nh^d} \sum_i Z_i, \qquad Z_i = K_i Y_i - \mathbb{E}K_i Y_i, K_i = K\left(\frac{X_i - x}{h}\right).$$

#### 4.1.1. Independence

**Proposition 6.** *Assuming that $(X_i, Y_i)$ is i.i.d. and that (A2), (A3) with $r = q$ hold, then the (A4) moment inequalities hold.*

*4.1.2. Non-causal weak dependence*

We now work as in [14], except for the necessary truncation used in [2].

**Proposition 7.** *Assume that $(X_i, Y_i)$ is $\lambda$-weakly dependent and that* (A2) *and* (A3) *hold. Assume that $\|Y_0\|_s < \infty$ for some $s > 2p$. If $\lambda(i) = \mathcal{O}(i^{-\lambda})$ with $\lambda > \frac{r(2r(s-p)+2p-s)}{(r-p)(s-2p)(r-1)}(p-1) \vee \frac{2(d-1)}{d}(q-1)$, then the* (A4) *moment inequalities hold.*

*4.1.3. Strong mixing*

**Proposition 8.** *Assume that* (A2) *and* (A3) *with $r > q$ hold, and that $(X_i, Y_i)$ is $\alpha$-mixing with $\alpha_i = \mathcal{O}(i^{-\alpha})$. If we also suppose that either:*

- $\alpha > ((q-1)\frac{r}{r-q}) \vee \frac{4sr-2s-4r}{(r-2)(s-4)}$ *and $h \sim n^{-a}$ with $ad \leq \frac{1-2/p}{3-2/r}$, or*
- $p, q$ *are even integers and $\alpha > \frac{r}{2}\frac{s-2p}{s-p}(1 - \frac{1}{p})$,*

*then the* (A4) *moment inequalities hold.*

**Remarks.**

- In the first item, the previous limitation on $h$ can be expressed as $p \geq d/\rho + 2$ if one makes use of the window width $h \sim n^{-1/(2\rho+d)}$, optimal with respect to power loss functions; this loss does not appear for integral order moments of the second item.
- In the case of absolute regularity, Viennet [29] provides sharp bounds for some integrals of the second order moments of such expressions. We do not derive them here, even if integrated square errors have specific interpretations: we need higher order moments in our case.

## 4.2. Uniform mean estimates

We now investigate uniform bounds:

$$\delta_n(B) = \left\| \sup_{x \in B} |\widehat{r}(x) - r(x)| \right\|_p. \tag{8}$$

In this setting, Ango Nze and Doukhan [1] prove the needed results under mixing assumptions; Ango Nze *et al.* [4] and Ango Nze and Doukhan [3] provide bounds under weak dependence conditions. For this, assumptions and lemmas need to be rephrased by replacing $\widehat{N}_n - N_n$ and $\widehat{D}_n - D_n$ by suprema of those expressions over $x \in B$ for some compact subset $B \subset \mathbb{R}^d$:

(A5) The condition (A1) holds for each $x \in B$.

(A6) The condition (A3) holds for each $x \in B$.

(A7) For some $q > p$ and $w_n = \frac{\sqrt{\log n}}{\sqrt{nh^d}} = v_n\sqrt{\log n}$, there exists $c > 0$ such that:

  1. $\|\sup_{x \in B} |\widehat{f}(x) - \mathbb{E}\widehat{f}(x)|\|_q \leq cw_n$;

2. $\|\sup_{x \in B} |\widehat{g}(x) - \mathbb{E}\widehat{g}(x)|\|_p \leq c w_n$.

We begin with two preliminary propositions before stating our main result.

**Proposition 9 (Uniform bias).** *Assuming that $f$ is bounded below over an open neighborhood of $B$ and that* (A2) *and* (A5) *hold, we have*

$$\sup_{x \in B} \left| r(x) - \frac{\mathbb{E}\widehat{g}(x)}{\mathbb{E}\widehat{f}(x)} \right| = \mathcal{O}(h_n^\rho).$$

**Proposition 10.** *Assuming that $f$ is bounded below over an open neighborhood of $B$ and that assumptions* (A2), (A6) *and* (A7) *hold, we have*

$$\left\| \sup_{x \in B} \left| \widehat{r}(x) - \frac{\mathbb{E}\widehat{g}(x)}{\mathbb{E}\widehat{f}(x)} \right| \right\|_p \leq C(1 + h^{d\beta(1/r-1)}(nh^d)^{(\beta-1)/2} + (nh^d)^{-\alpha/2} n^{1/s}) w_n$$

*with $\alpha = q(\frac{1}{p} - \frac{1}{s} - \frac{1}{q})$ and $\beta = \frac{pr}{q(r-p)}$.*

The following theorem derives from the two previous propositions.

**Theorem 3.** *Let $(X_t, Y_t)$ be a stationary sequence. Assume that conditions* (A2), (A5), (A6) *and* (A7) *hold for some $s > 2p$. The optimal rate for $h$ is*

$$h = C \left( \frac{\log n}{n} \right)^{1/(2\rho+d)} \tag{9}$$

*and*

$$\left\| \sup_{x \in B} |\widehat{r}(x) - r(x)| \right\|_p \leq C \left( \frac{\log n}{n} \right)^{\rho/(2\rho+d)}. \tag{10}$$

The end of the section is devoted to different dependence conditions that are sufficient for (A7) to hold.

### 4.2.1. Independence

First, we evaluate the uniform bound for the moments under independence.

**Proposition 11.** *Let $(X_t, Y_t)_{t \in \mathbb{N}}$ be an i.i.d. sequence. If we assume that conditions* (A2), (A5) *and* (A6) *hold for some $s > 2p$, $\rho > dp/(s-2p)$, then* (A7) *holds, hence* (9) *yields the bounds*

$$\left\| \sup_{x \in B} |\widehat{f}(x) - f(x)| \right\|_p \leq C (\log n/n)^{\rho/(2\rho+d)}, \tag{11}$$

$$\left\| \sup_{x \in B} |\widehat{g}(x) - g(x)| \right\|_p \leq C (\log n/n)^{\rho/(2\rho+d)}. \tag{12}$$

*4.2.2. Absolute regularity*

**Proposition 12.** *Let $(X_i, Y_i)_{i \in \mathbb{N}}$ be an absolute regular (also called $\beta$-mixing) sequence. Assume that conditions* (A2), (A5) *and* (A6) *hold for some $s > 2p$, $\rho > dp/(s-2p)$. If we assume that the mixing coefficients satisfy $\beta_i = \mathcal{O}(i^{-\beta})$ with $\beta > \frac{s\rho + (2s-p)d}{\rho(s-2p)-pd} \vee (1 + \frac{2d}{\rho})$, then assumption* (A7) *holds, hence the choice (9) yields the bounds (11)–(12).*

*4.2.3. Strong mixing*

Using the Fuk–Nagaev inequality in [24] also yields an analogous result.

**Proposition 13.** *Assume that the process $(X_i, Y_i)_{i \in \mathbb{N}}$ is stationary and strongly mixing with $\alpha_i = \mathcal{O}(i^{-\alpha})$ for $\alpha > \frac{3\rho s + 2ds + d\rho s - 4\rho p - 3dp - d\rho p}{dp - \rho(s-2p)} \vee 2\frac{s-1}{s-2}$. If we further assume that conditions* (A2), (A5) *and* (A6) *hold for some $s > 2p$, $\rho > dp/(s-2p)$, then assumption* (A7) *holds, hence (9) yields the bounds (11)–(12).*

*4.2.4. Non-causal weak dependence*

**Proposition 14.** *Assume that the process $(X_i, Y_i)_{i \in \mathbb{Z}}$ is stationary and $\lambda$-weakly dependent with $\lambda(i) = \mathcal{O}(e^{-\lambda i^b})$, $b > 0$. If we further assume that conditions* (A2), (A5) *and* (A6) *hold for some $s > 2p$, $\rho > dp/(s-2p)$, then assumption* (A7) *holds, hence (9) yields the bounds (11)–(12).*

*Remark.* Other dependence settings may also be addressed. For example, the $\phi$-mixing case considered in [8] and the use of coupling in weakly dependent sequences by Dedecker and Prieur [10] both yield suitable exponential inequalities to complete analogous results.

## 5. Proofs

In the proofs, $C > 0$ is a constant which may change from one line to another.

### 5.1. Proof of the main lemmas

**Proof of Lemma 2.** Setting $z = (D_n - \widehat{D}_n)/D_n$, we rewrite

$$\Delta_n = \frac{\widehat{N}_n}{D_n} \cdot \frac{1}{1-z} - \frac{N_n}{D_n} = \frac{\widehat{N}_n - N_n}{D_n} + \frac{\widehat{N}_n}{D_n}\left(\frac{1}{1-z} - 1\right). \tag{13}$$

But $\frac{1}{1-z} - 1 = \frac{z}{1-z} = z + \frac{z^2}{1-z}$, hence for $\alpha \in [0,1]$,

$$\left|\frac{1}{1-z} - 1\right| \leq |z| + \frac{|z| \wedge |z|^2}{|1-z|} \leq |z| + \frac{|z|^{1+\alpha}}{|1-z|},$$

which implies that

$$D_n|\Delta_n| \leq |\widehat{N}_n - N_n| + |\widehat{N}_n|\frac{|\widehat{D}_n - D_n|}{D_n} + |\widehat{N}_n|\left|\frac{D_n}{\widehat{D}_n}\right|\frac{|\widehat{D}_n - D_n|^{1+\alpha}}{|D_n|^{1+\alpha}}$$

$$\leq |\widehat{N}_n - N_n| + \frac{|\widehat{N}_n|}{D_n}|\widehat{D}_n - D_n| + \left|\frac{\widehat{N}_n}{\widehat{D}_n}\right|\frac{|\widehat{D}_n - D_n|^{1+\alpha}}{|D_n|^\alpha}$$

$$\leq |\widehat{N}_n - N_n| + \frac{|\widehat{N}_n|}{D_n}|\widehat{D}_n - D_n| + \max_{1 \leq i \leq n}|V_{i,n}|\frac{|\widehat{D}_n - D_n|^{1+\alpha}}{|D_n|^\alpha}.$$

**Proof of Lemma 1.** From the preceding relation, we have

$$D_n\|\Delta_n\|_p \leq v_n + \frac{1}{D_n}\|\widehat{N}_n(\widehat{D}_n - D_n)\|_p + \left\|\max_{1 \leq i \leq n}|V_{i,n}|\frac{|\widehat{D}_n - D_n|^{1+\alpha}}{|D_n|^\alpha}\right\|_p. \qquad (14)$$

From the Hölder inequality with exponents $1/a + 1/b = 1$, we have

$$\left\|\frac{|\widehat{N}_n|}{D_n}|\widehat{D}_n - D_n|\right\|_p \leq \frac{1}{D_n}\|\widehat{N}_n\|_{pa}\|\widehat{D}_n - D_n\|_{pb}.$$

Now, the assumption $\|U_{i,n}V_{i,n}\|_r \leq C_n$ implies that $\|\widehat{N}_n\|_r \leq C_n$. The second term in the right-hand side of inequality (14) is bounded using the property $\|\widehat{N}_n\|_{pa} \leq |N_n| + \|\widehat{N}_n - N_n\|_{pa}$. Consider now some $\beta \in [0,1]$ and $u, v \geq 0$ such that $1/u + 1/v = 1$, to be determined later. Then $|\widehat{N}_n - N_n| = |\widehat{N}_n - N_n|^\beta|\widehat{N}_n - N_n|^{1-\beta}$; the Hölder inequality implies that if we choose $upa(1 - \beta) = p$ and $vpa\beta = r$, then

$$\|\widehat{N}_n - N_n\|_{pa} \leq \|\widehat{N}_n - N_n\|_{upa(1-\beta)}^{1-\beta}\|\widehat{N}_n - N_n\|_{vpa\beta}^\beta$$

$$\leq \|\widehat{N}_n - N_n\|_p^{1-\beta}\|\widehat{N}_n - N_n\|_r^\beta$$

$$\leq v_n^{1-\beta}(|N_n|^\beta + C_n^\beta),$$

thus $\|\widehat{N}_n\|_{pa} \leq |N_n| + v_n^{1-\beta}(|N_n|^\beta + C_n^\beta)$. Setting $b = q/p$, we derive $a = (q-p)/q$, hence $\frac{1}{u} = \frac{q}{q-p}(1-\beta)$ and $\frac{1}{v} = \frac{pq}{r(q-p)}\beta$. With the relation $1/u + 1/v = 1$, we find $\beta = pr/q(r-p)$. Then

$$\frac{\|\widehat{N}_n\|_{pa}}{D_n}\|\widehat{D}_n - D_n\|_{pb} \leq \frac{1}{D_n}(|N_n| + v_n^{1-\beta}(|N_n|^\beta + C_n^\beta))v_n.$$

The last term in relation (14) is more difficult to handle; it may be bounded using the Hölder inequality with exponents $1/a + 1/b = 1$ and

$$\left\|\max_{1 \leq i \leq n}|V_{i,n}|\frac{|\widehat{D}_n - D_n|^{1+\alpha}}{|D_n|^\alpha}\right\|_p \leq \frac{1}{|D_n|^\alpha}\left\|\max_{1 \leq i \leq n}|V_{i,n}|\right\|_{pa}\||\widehat{D}_n - D_n|^{1+\alpha}\|_{pb}$$

$$\leq \frac{1}{|D_n|^\alpha}\Big(\mathbb{E}\max_{1\leq i\leq n}|V_{i,n}|^{pa}\Big)^{1/(pa)}v_n^{1+\alpha},$$

if $q \geq pb(1+\alpha)$ or, equivalently, if $q - p(1+\alpha) \geq q/a$. We use an argument of Pisier [22]: if $\varphi:\mathbb{R}^+ \to \mathbb{R}^+$ is convex and non-decreasing, then

$$\varphi\Big(\mathbb{E}\max_i |V_{i,n}|^{pa}\Big) \leq \mathbb{E}\varphi\Big(\max_i |V_{i,n}|^{pa}\Big) \tag{15}$$
$$\leq \mathbb{E}\sum_i \varphi(|V_{i,n}|^{pa}) \leq \sum_i \mathbb{E}\varphi(|V_{i,n}|^{pa}).$$

Hence $\mathbb{E}\max_i |V_{i,n}|^{pa} \leq (nc)^{pa/s}$ with $\varphi(x) = x^{s/pa}$. Now, the bound in the right-hand side of (14) can be expressed as $v_n^{1+\alpha}(nc)^{1/s}/|D_n|^\alpha$ if $s \geq pa$; for this, we use $1 - \frac{p}{q}(1+\alpha) \geq \frac{1}{a} \geq \frac{p}{s}$ if $\alpha > 0$ is small enough with $\frac{1}{p} \geq \frac{1+\alpha}{q} + \frac{1}{s}$. □

**Proof of Corollary 1.** The proof of this result is standard. Namely, $p' < p$ implies from the Markov inequality that for $Z_n = v_n^{-1}|\Delta_n|$, the sequence $(Z_n)^{p'}$ is uniformly integrable. Set $f_k(z) = |z|^{p'}\mathbb{1}_{\{|z|\leq k\}}$. From convergence in distribution,

$$\mathbb{E}f_k(Z_n) \underset{n\to\infty}{\to} \mathbb{E}f_k(Z)$$

and this occurs uniformly with respect to $k$ from uniform integrability. □

### 5.2. Proofs for Section 3

**Proof of Theorem 1.** Set $\Delta_n = \widehat{R}_n - R$ and refer to Lemma 1: here $r$ and $s$ are such that $\beta = 1$, and $\alpha = 2/s$. Because $v_n = Cn^{-1/2}$, we get

$$D\|\Delta_n\|_p \leq \left(1 + 2\frac{N}{D} + \frac{c}{D} + \frac{cC^\alpha n^{-\alpha/2+1/s}}{D^\alpha}\right)v_n,$$

where the last term in the parenthesis is bounded with respect to $n$, implying that

$$\|\widehat{R}_n - R\|_p = \mathcal{O}(1/\sqrt{n}).\qquad\square$$

**Proof of Proposition 1.** Set $Z_i = U_iV_i - \mathbb{E}U_iV_i$. Let $\alpha^{-1}(u) = \sum_{i\geq 0}\mathbb{1}\{u < \alpha_i\}$, and denote by $Q_Z$ the generalized inverse of the tail function $x \mapsto \mathbb{P}(|Z_0| > x)$. From heredity, the mixing coefficient $\alpha_{Z_i}$ of the sequence $Z$ is bounded by $\alpha_i$. Theorem 2.5 in [24] shows that

$$\left\|\sum_{i=1}^n Z_i\right\|_p \leq \sqrt{2pn}\left(\int_0^1 (\alpha^{-1}(u)\wedge n)^{p/2}Q_Z^p(u)\,\mathrm{d}u\right)^{1/p}.$$

But, as $\|Z_0\|_r \leq c$, we have $\mathbb{E}|Z_0|^r = \int_0^\infty Q_Z^r(u)\,du \leq c^r$. Define $a = r/p$ and $b = r/(r-p)$. Then, by the Hölder inequality,

$$\left\|\sum_{i=1}^n Z_i\right\|_p \leq \sqrt{2pn}\left(\int_0^1 (\alpha^{-1}(u) \wedge n)^{pb/2}\,du\right)^{1/(bp)}\left(\int_0^1 Q_Z^r(u)\,du\right)^{1/(ap)}.$$

From $\alpha > pb/2$, the first integral is finite and $\|\widehat{N}_n - N\|_p \leq C/\sqrt{n}$. For $Y_n = U_n - \mathbb{E}U_n$ and $b' = r'/(r'-q)$, we similarly get that $\|\widehat{D}_n - D\|_q \leq C/\sqrt{n}$. □

**Proof of Proposition 2.** If $U$ and $V$ are independent sequences, then $(Y_i) = (U_i - \mathbb{E}U_i)$ is $\gamma$-dependent with the same coefficients, as is $(Z_i) = (U_i - \mathbb{E}U_i)(V_i - \mathbb{E}V_i)$, since $\|\mathbb{E}(Z_i|\mathcal{M}_0)\|_1 \leq \|\mathbb{E}(U_i|\mathcal{M}_0) - \mathbb{E}U_i\|_1 \|\mathbb{E}(V_i|\mathcal{M}_0) - \mathbb{E}V_i\|_1 \leq 2c\gamma_i$. Denote by $G_Z$ the inverse of $x \mapsto \int_0^x Q_Z(u)\,du$. Corollary 5.3, page 124 in [9] states that

$$\left\|\sum_{i=1}^n Z_i\right\|_p \leq \sqrt{2pn}\left(\int_0^{\|Z\|_1} (\gamma^{-1}(u) \wedge n)^{p/2} Q_Z^{p-1} \circ G_Z(u)\,du\right)^{1/p}.$$

As $U$ is bounded by $c$, $\mathbb{E}|Z|^s = \int_0^\infty Q_Z^{s-1} \circ G_Z(u)\,du \leq c^s$.

Define $a = (s-1)/(p-1)$ and $b = (s-1)/(s-p)$. Then, by the Hölder inequality,

$$\left\|\sum_{i=1}^n Z_i\right\|_p \leq \sqrt{2pn}\left(\int_0^\infty (\gamma^{-1}(u) \wedge n)^{pb/2}\,du\right)^{1/(bp)}\left(\int_0^\infty Q_Z^{r-1} \circ G_Z(u)\,du\right)^{1/(ap)}.$$

Because $\gamma > \frac{pb}{2}$, the first integral is finite and $\|\widehat{N}_n - N\|_p \leq \frac{C}{\sqrt{n}}$. Similarly,

$$\left\|\sum_{i=1}^n Y_i\right\|_q \leq \sqrt{2qn}\left(\int_0^{\|Y\|_1} (\gamma^{-1}(u) \wedge n)^{q/2} Q_Y^{q-1} \circ G_Y(u)\,du\right)^{1/q}$$

$$\leq \sqrt{2qn}(2c)^{q-1}\left(\int_0^{2c} (\gamma^{-1}(u) \wedge n)^{q/2}\,du\right)^{1/q}.$$

Because $\gamma > \frac{q}{2}$, the first integral is finite and $\|\widehat{D}_n - D\|_q \leq \frac{C}{\sqrt{n}}$. □

**Proof of Proposition 3.** Let $Z_i = U_i V_i - \mathbb{E}U_i V_i$.

- Define $Z_I = (Z_{i_1}, \ldots, Z_{i_u})$. If $U$ and $V$ are independent, then

$$\operatorname{Cov}(g_1(Z_I), g_2(Z_J)) = \operatorname{Cov}(\mathbb{E}(g_1(Z_I)|U), \mathbb{E}(g_2(Z_J)|U))$$
$$+ \mathbb{E}\operatorname{Cov}(g_1(Z_I), g_2(Z_J)|U).$$

If we define $\tilde{g}_1(u_I) = \mathbb{E}(g_1(Z_I)|U_I = u_I)$, then

$$|\tilde{g}_1(u_I) - \tilde{g}_1(u_I')| \leq \operatorname{Lip} g_1 \mathbb{E}|V_0|\sum_{i \in I}|u_i - u_i'|$$

so that $g_1(U_I)$ is Lipschitz with respect to $(U_I, V_I)$. We have

$$\text{Cov}(g_1(Z_I), g_2(Z_J)|U_{I\cup J} = u_{I\cup J}) = \text{Cov}(\check{g}_1(V_I), \check{g}_2(V_J)).$$

Here, $\check{g}_1$ is a Lipschitz function with coefficient $\max_{i\in I} u_i \operatorname{Lip} g_1$. From (6), we conclude that $(Z_i)$ is also $\lambda$-weakly dependent with $\lambda_Z(k) \leq (c \wedge \mathbb{E}V_0)k^{-\lambda}$. From Theorem 4.1, page 77, and Proposition 13.1, page 293 in [9] extended to $\lambda$-dependence,

$$\left\| \sum_{i=1}^n Z_i \right\|_q^q \leq \frac{C(2q-2)!}{(q-1)!} \left\{ \left( n \sum_{k=0}^{n-1} \lambda_Z(k) \right)^{q/2} \vee \left( M^{q-2} n \sum_{k=0}^{n-1} (k+1)^{q-2} \lambda_Z(k) \right) \right\}.$$

The second term is negligible as $n$ tends to infinity and the first sum over $k$ is bounded so that $\|\sum_{i=1}^n Z_i\|_q \leq Cn^{1/2}$ and thus $\|\widehat{N}_n - N\|_q \leq \frac{C}{\sqrt{n}}$.

- In the second case, $\|U_0\|_{r'} + \|V_0\|_{r'} \leq \infty$ and Proposition 2.1, page 12 in [9] implies that $\lambda_Z(k) \leq Ck^{-((r'-2)\lambda)/r'}$. □

### 5.3. Proofs for Section 4.1

**Proof of Lemma 4.** The previous convergences $\mathbb{E}\widehat{f}(x) \to f(x)$ and $\mathbb{E}\widehat{g}(x) \to r(x)f(x)$ are controlled by $\mathcal{O}(h_n^\rho)$ under $\rho$-regularity conditions (A1) (see Ango Nze and Doukhan [3]). Write

$$\left| r(x) - \frac{\mathbb{E}\widehat{g}(x)}{\mathbb{E}\widehat{f}(x)} \right| \leq \frac{|g(x) - \mathbb{E}\widehat{g}(x)|}{f(x)} + |\mathbb{E}\widehat{g}(x)| \frac{|\mathbb{E}\widehat{f}(x) - f(x)|}{f(x)\mathbb{E}\widehat{f}(x)}.$$

Since $f$ is bounded below by 0 around $x$, $\mathbb{E}\widehat{f}(x)$ is also bounded below by 0 and we get the result. □

**Proof of Proposition 5.** Condition ((A3)-1) gives $\|V_{i,n}\|_s < \infty$ and $\|U_{i,n}V_{i,n}\|_r = \mathcal{O}(h^{d(1/r-1)})$ for $r < s$ whenever ((A3)-3) holds.

Indeed, $\|V_{i,n}\|_s = \|Y_0\|_s = c < \infty$, so, for $r \leq s$,

$$\mathbb{E}|U_{i,n}V_{i,n}|^r \leq \int\int h^{-rd} K^r((X_i - x)/h) |y|^r f(x,y) \, dx \, dy$$

$$\leq \int h^{-rd} K^r((X_i - x)/h) g_r(x) \, dx$$

$$\leq \|g_r\|_\infty \int h^{-rd} K^r((X_i - x)/h) \, dx$$

$$\leq Ch^{(1-r)d}.$$

Set $K_{h_n}(\cdot) = h_n^{-d} K(\cdot/h_n)$ and denote by $\star$ the convolution. Here, $D_n = \mathbb{E}\widehat{f}(x) = \mathbb{E}U_{i,n} = f \star K_{h_n}(x) \to f(x)$, the marginal density of $X_0$, and $N_n = \mathbb{E}\widehat{g}(x) = \mathbb{E}U_{i,n}V_{i,n} = (rf) \star$

$K_{h_n}(t) \to r(x)f(x)$. From Lemma 1, we get

$$D_n \left\| \widehat{r}(x) - \frac{\mathbb{E}\widehat{g}(x)}{\mathbb{E}\widehat{f}(x)} \right\|_p \leq \left(1 + \frac{|N_n|}{D_n} + \frac{|N_n|^\beta v_n^{1-\beta}}{D_n} + \frac{C_n^\beta v_n^{1-\beta}}{D_n} + \frac{v_n^\alpha c_n n^{1/s}}{D_n^\alpha}\right) v_n,$$

where $D_n$, $N_n$ and $c_n$ are equivalent to constants, $C_n \equiv Ch^{d(1/r-1)}$ and $v_n \equiv C(nh^d)^{-1/2}$. Substituting the orders in the different terms, we obtain the result. □

**Proof of Theorem 2.** From the preceding propositions, we get

$$\|\widehat{r}(x) - r(x)\|_p \leq Ch_n^\rho + C(1 + h^{d\beta(1/r-1)}(nh^d)^{(\beta-1)/2} + (nh^d)^{-\alpha/2}n^{1/s})v_n. \tag{16}$$

Note that $h_n^\rho = Cv_n$. The expression in parentheses is bounded if

$$0 \leq \frac{\beta d(1-r)}{r} + (1-\beta)\rho, \qquad 0 \leq \frac{\alpha\rho}{d+2\rho} - \frac{1}{s}.$$

These conditions correspond to (7). □

*A bound of interest which does not use dependence conditions*

The proofs of the propositions under different kinds of dependence make use of a common bound that holds in all cases. For a positive integer $k$, we define the coefficients of weak dependence as non-decreasing sequences $(C_{k,q})_{q \geq 2}$ such that

$$C_{k,q} = \sup |\operatorname{Cov}(Z_{i_1} \cdots Z_{i_m}, Z_{i_{m+1}} \cdots Z_{i_q})|, \tag{17}$$

where the supremum is taken over all $\{i_1, \ldots, i_q\}$ such that $1 \leq i_1 \leq \cdots \leq i_q$ and where $m$, $k$ satisfy $i_{m+1} - i_m = k$. Independently of the dependence structure, we get a bound for $C_{k,q}$.

**Lemma 3.** *If we assume* (A3) *and* (A4), *then* $C_{k,p} \leq Ch^{2d(s-p)/s}$.

**Proof.** Define $\{i_1, \ldots, i_p\}$ as a sequence that attains the sup defining $C_{k,p}$.

$$\begin{aligned}
C_{k,p} &= |\mathbb{E}Z_{i_1} \cdots Z_{i_p}| \\
&\leq 2^p \mathbb{E}|Y_{i_1} K_{i_1} \cdots Y_{i_p} K_{i_p}| \\
&\leq 2^p (\mathbb{E}|Y_{i_1} \cdots Y_{i_p}|^{s/p})^{p/s} \left(\max_{i_1, i_2} \mathbb{E}(K_{i_1} K_{i_2})^{s/(s-p)}\right)^{1-p/s} \\
&\leq 2^p \|Y_0\|_s^p \left(\int\int K^{s/(s-p)}\left(\frac{x-u}{h}\right) K^{s/(s-p)}\left(\frac{x-t}{h}\right) G(u,t) \, du \, dt\right)^{1-p/s} \\
&\leq 2^p \|Y_0\|_s^p h^{2d(1-p/s)}.
\end{aligned}$$

The claim in the remark following (A3) is based on the fact that there is no need to use the Hölder inequality if $H$ is bounded around $(x, x)$ and thus $C_{k,p} \leq Ch^{2d}$. Now $C_{k,q} \leq Ch^{2d}$ also holds for the denominator (we may also set $Y_i \equiv 1$). □

**Proof of Proposition 6.** From the Rosenthal inequality for independent variables, there exist constants $C_q \geq 0$ depending only on $q$ (see, for example, Figiel *et al.* [19] for more details concerning the constants) such that

$$\mathbb{E}\left|\sum_{i=1}^n Z_i\right|^q \leq C_q\left(\left|\sum_{i=1}^n \mathbb{E}Z_i^2\right|^{q/2} + \sum_{i=1}^n \mathbb{E}|Z_i|^q\right).$$

Here, $\mathbb{E}Z_i^q \leq 2^q \mathbb{E}|K_i Y_i|^q$. In the beginning of the proof of Proposition 5, we get $\mathbb{E}|K_i Y_i|^q \leq \|g_q\|_\infty \|K\|_q^q h^d$, and we deduce $\mathbb{E}|\sum Z_i|^q \leq C((nh^d)^{q/2} + nh^d)$, and

$$\|\widehat{g}(x) - \mathbb{E}\widehat{g}(x)\|_q \leq \frac{C}{\sqrt{nh^d}}.$$

The case of the denominator is obtained by setting $Y_i \equiv 1$. □

**Proof of Proposition 7.** We first establish a Rosenthal inequality for weakly dependent variables. For any integer $p \geq 2$, $\mathbb{E}(\sum_{i=0}^n Z_i)^p \leq p! A_p$, where

$$A_p = \sum_{1 \leq i_1 \leq \cdots \leq i_p \leq n} |\mathbb{E}(Z_{i_1} \cdots Z_{i_p})|. \tag{18}$$

Then

$$A_p \leq n \sum_{k=0}^n (k+1)^{p-2} C_{k,p} + \sum_{l=2}^{p-2} A_l A_{p-l}. \tag{19}$$

To bound the sum with coefficients $C_{k,p}$, we shall use the following lemma.

**Lemma 4 (Doukhan and Neumann [15], Lemma 10-(11)).** *If we assume that the stationary sequence $(Z_n)_{n\in\mathbb{Z}}$ is $\lambda$-weakly dependent and satisfies $\mu = \|Z_0\|_r \leq 1$ for some $r > p$, then*

$$C_{k,p} \leq 2^{p+3} p^4 \mu^{r(p-1)/(r-1)} \lambda(k)^{(r-p)/(r-2)}.$$

Recall that if $(X_i, Y_i)$ is $\lambda$-weakly dependent with $\lambda(k) \leq Ck^{-\lambda}$, then $Z_i$ is $\lambda$-weakly dependent with $\lambda(k) \leq Ck^{-\lambda(r-2)/r}$ (see Proposition 2.1, page 12 in our monograph [9]). From Lemma 4, we get $C_{k,p} \leq Ch^{d(p-1)/(r-1)}k^{-\lambda(1(p)/r)}$.

We define conditions on the dependence coefficients that ensure that the sum to control is of the same order of magnitude as its first term, that is, $\mathcal{O}(h^d)$. From Lemma 3, for any $0 < b < 1$,

$$C_{k,p} \leq Ch^d(h^{-d(r-p)/(r-1)}k^{-\lambda(r-p)/r} \wedge h^{d(1-2p/s)})$$

$$\leq Ch^d \times k^{-\lambda b(r-p)/r} h^{d(-b(r-p)/(r-1)+(1-b)(1-2p/s))}.$$

Choosing $b < \frac{(s-2p)(r-1)}{2r(s-p)+2p-s}$, the exponent of $h$ in the parentheses is positive. Because $\lambda > \frac{1}{b} \cdot \frac{r(p-1)}{r-p}$, the sum over $k$ (in inequality (19)) converges and is less than a constant, say $a_p$. We get $A_p \leq a_p n h^d + \sum_{l=2}^{p-2} A_l A_{p-l}$ and from $A_2 \leq a_2 n h^d$, we deduce by induction on $p$ that $A_p \leq c_p (nh^d)^{p/2}$ for a sequence $c_p$ which may also depend on $x$. Hence, $\|\widehat{g}(x) - \mathbb{E}\widehat{g}(x)\|_p \leq c_p(x)(nh^d)^{-p/2}$.

The case of the denominator is obtained by setting $Y_i \equiv 1$. In this case, we also note that the bound in Lemma 3 can be expressed as $C_{k,q} \leq Ch^{2d}$, since a Hölder inequality is no longer needed.

$$C_{k,q} \leq Ch^d(h^{-2+d}k^{-\lambda} \wedge h^d)$$
$$\leq Ch^d \times k^{-\lambda b} h^{(-b(d-2)+(1-b)d)}.$$

Choosing $b < d/(2d-2)$, the exponent of $h$ in the parentheses is positive. Because $\lambda > q - 1/b$, the sum over $k$ in inequality (19) converges. $\square$

**Proof of Proposition 8.** • Under the first set of conditions, Rio [24], Theorem 6.3 states the following Rosenthal inequality:

$$\mathbb{E}\Big|\sum Z_i\Big|^p \leq a_p \bigg(\sum_i \sum_j |\operatorname{Cov}(Z_i, Z_j)|\bigg)^{p/2} + nb_p \int_0^1 (\alpha^{-1}(u) \wedge n)^{p-1} Q_Z^p(u)\,du.$$

We use Lemma 3 to prove that the first term is $\mathcal{O}((nh^d)^{p/2})$. From the Davydov inequality, we get a second bound for the covariance:

$$|\operatorname{Cov}(Z_0, Z_i)| \leq 6\alpha_i^{(r-2)/r} \|Y_0 K_0\|_r^2 \leq C\alpha_i^{(r-2)/r} h^{2d(1/r-1)},$$

$$\operatorname{Var}(Z_0) \sim h^d g_2(x) \int K^2(u)\,du.$$

Hence,

$$\bigg|\sum_i \sum_j \operatorname{Cov}(Z_i, Z_j) - \sum_i \operatorname{Var}(Z_i)\bigg| \leq nh^d \sum_i \alpha_i^{(r-2)/r} h^{d(2/r-3)} \wedge h^{d(1-4/s)}.$$

Thus, considering some $0 < b < 1$,

$$\sum_i \alpha_i^{1-2/r} h^{d(2/r-3)} \wedge h^{d(1-4/s)} \leq h^{d(2b(1/r+2/s-2)+1-4/s)} \sum_i \alpha_i^{b(1-2/r)}.$$

This last term tends to 0 if $b > \frac{sr-4r}{4sr-2s-4r}$ and $\alpha > \frac{r}{b(r-2)}$. This is possible if $\alpha > \frac{r(4sr-2s-4r)}{(r-2)(sr-4r)}$.

Consider the second term and apply the Hölder inequality with exponents $r/(r-p)$ and $r/p$:

$$n \int_0^1 (\alpha^{-1}(u) \wedge n)^{p-1} Q_Z^p(u) \, du \leq n \left( \int_0^1 (\alpha^{-1}(u))^{r(p-1)/(r-p)} \, du \right)^{(r-p)/r} \|Z\|_r^p.$$

The first integral is convergent as soon as $\alpha > \frac{r(p-1)}{r-p}$ and, from assumption (A3)-3, $\|Z\|_r^r \leq Ch^{d(1-r)}$ so that this second term is negligible if $nh^{dp(1-r)/r} \leq C(nh^d)^{p/2}$. Hence, if the sequence $n^{r(2-p)} h^{dp(2-3r)}$ is bounded, we obtain the desired bound.

Now, consider the denominator,

$$\left| \sum_i \sum_j \mathrm{Cov}(K_i, K_j) - \sum_i \mathrm{Var}(K_i) \right| \leq nh^d \sum_i \alpha_i h^{-d} \wedge h^d.$$

Thus, $\sum_i \alpha_i h^{-d} \wedge h^d \leq h^{d(1-2b)} \sum_i \alpha_i^b$ for $0 < b < 1$. This last term tends to 0 if $b > 1/2$, which implies that $\alpha > 2$.

Consider the second term,

$$n \int_0^1 (\alpha^{-1}(u) \wedge n)^{q-1} Q_{K_0}^q(u) \, du \leq n \left( \int_0^1 (\alpha^{-1}(u))^{r(q-1)/(r-q)} \, du \right)^{(r-q)/r} \|K_0\|_r^q.$$

The first integral is finite if $\alpha > \frac{r(q-1)}{r-q}$. Analogously, the second term is negligible as soon as $n^{r(2-q)} h^{dq(2-r)}$ is a bounded sequence.

Hence, if $h \sim n^{-a}$, a monotonicity argument shows that the previous bounds require $ad \leq \frac{1-2/p}{3-2/r}$.

• Under the second set of conditions, we use the idea from the proof of Proposition 7 (this idea was initiated in [16]). We again use relations (19), and expression (18) is bounded by using the alternative bound of $C_{k,p}$, which can be expressed as $|\mathbb{E} Z_{i_1} \cdots Z_{i_p}|$ for a suitable sequence $i_1 \leq \cdots \leq i_p$ with $i_{u+1} - i_u = k$. The Davydov inequality (see Theorem 3(i) in [13]) and the Hölder inequality then together imply that

$$C_{k,p} \leq 6 \alpha_k^{p/r} \|Z_{i_1} \cdots Z_{i_u}\|_{r/u} \|Z_{i_{u+1}} \cdots Z_{i_p}\|_{r/(p-u)} \leq 6 \alpha_k^{p/r} \|Z_0\|_r^p$$

and $C_{k,p} \leq Ch^d (\alpha_k^{p/r} h^{-d}) \wedge h^{d(1-2p/s)} \leq Ch^d \alpha_k^{bp/r} h^{d(1-2p/s-2b(1-p/(2s)))}$ from Lemma 3 if $0 \leq b \leq 1$. Then, setting $b = (s-2p)/(2(s-p))$,

$$n \sum_{k=0}^n (k+1)^{p-2} C_{k,p} = \mathcal{O}(nh^d) \qquad \text{if } \alpha > \frac{r}{2} \frac{s-2p}{s-p} \left(1 - \frac{1}{p}\right).$$

The case of the denominator is exactly analogous and here we replace $=p$ by $q$, $s$ by $\infty$ and, in order to let the previous condition unchanged, we replace $r$ by $r'$ with $\frac{r'}{2}(1 - \frac{1}{q}) = \frac{r}{2}\frac{s-2p}{s-p}(1 - \frac{1}{p})$. □

## 5.4. Proofs for Section 4.2

**Proof of Proposition 9.** The previous convergences $\mathbb{E}\widehat{f}(x) \to f(x)$ and $\mathbb{E}\widehat{g}(x) \to r(x)f(x)$ are uniformly controlled by $\mathcal{O}(h_n^\rho)$ under $\rho$-regularity conditions (A5) from the continuity of derivatives over the considered sets and a standard compactness argument. The proofs do not use the positivity of $K$ so that arbitrary values for $\rho$ are possible; see Ango Nze and Doukhan [3]. Note that if $V \ni x$ denotes an open set over which the previous assumptions (A1) hold and such that $\inf_B f > 0$, then for each open set $W$ with $\overline{W} \subset V$, the previous relation holds uniformly over $W$. Hence, under (A5), if $V$ denotes an open set with $B \subset V$ such that the assumptions (A1) still hold, the bounds for biases hold uniformly over $B$. We thus proceed as in Proposition 4 to complete the proof. $\square$

**Proof of Proposition 10.** From Lemma 2,

$$\inf_{x \in B} D_n(x) \sup_{x \in B} \left| \widehat{r}(x) - \frac{\mathbb{E}\widehat{g}(x)}{\mathbb{E}\widehat{f}(x)} \right| \le \sup_B |\widehat{N}_n - N_n| + \sup_B \frac{|\widehat{N}_n|}{D_n} \sup_B |\widehat{D}_n - D_n|$$

$$+ \max_{1 \le i \le n} |Y_i| \sup_B \frac{|\widehat{D}_n - D_n|^{1+\alpha}}{|D_n|^\alpha}.$$

As in the proof of Lemma 1, substituting the supremum to the variables, we get

$$\inf_{x \in B} D_n(x) \left\| \sup_{x \in B} \left| \widehat{r}(x) - \frac{\mathbb{E}\widehat{g}(x)}{\mathbb{E}\widehat{f}(x)} \right| \right\|_p$$

$$\le C \left( 1 + \frac{\sup_B |N_n|}{\inf_B D_n} + \frac{\sup_B |N_n|^\beta w_n^{1-\beta}}{\inf_B D_n} + \frac{C_n^\beta w_n^{1-\beta}}{\inf_B D_n} + \frac{w_n^\alpha c_n n^{1/s}}{\inf_B D_n^\alpha} \right) w_n,$$

where $D_n$, $N_n$ and $c_n$ are equivalent to constants and $C_n \equiv Ch^{d(1/r-1)}$. Substituting the orders to the expressions gives the result. $\square$

*Truncation and variance estimation*

This paragraph introduces some common elements of proofs of Propositions 11–14.

Let $M > 0$ and consider the truncated modification of $Y_i$, $\widetilde{Y}_i = Y_i \mathbb{1}\{|Y_i| \le M\} - M\mathbb{1}\{Y_i < -M\} + M\mathbb{1}\{Y_i > M\}$. Define $\widetilde{g}(x) = \frac{1}{nh^d} \sum_{i=1}^n \widetilde{Y}_i K_i$. Then, from the Markov inequality,

$$\left\| \sup_{x \in \mathbb{R}^d} |\widehat{g}(x) - \widetilde{g}(x)| \right\|_p^p \le \frac{1}{h^{dp}} \mathbb{E}|Y|^p \mathbb{1}\{|Y| > M\} \le \frac{M^{p-s}}{h^{dp}} \mathbb{E}|Y|^s.$$

With the choice (9), in order to conveniently bound this term, we assume that

$$M \ge Ch^{-p/(s-p)(\rho+d)}. \tag{20}$$

Let $Z(x) = \sqrt{nh^d}|\widetilde{g}(x) - \mathbb{E}\widetilde{g}(x)|$. Below, we will need (uniform) bounds of

$$\operatorname{Var} Z(x) = \frac{1}{nh^d} \sum_{|i|<n} (n-|i|) \operatorname{Cov}(\widetilde{Y}_0 K_0, \widetilde{Y}_i K_i) \leq \frac{2}{h^d} \sum_{i=0}^{n-1} |\operatorname{Cov}(\widetilde{Y}_0 K_0, \widetilde{Y}_i K_i)|.$$

Set $\Gamma_i(x) = |\operatorname{Cov}(\widetilde{Y}_0 K_0, \widetilde{Y}_i K_i)|$. A first bound of $\Gamma_i(x)$ for $i \neq 0$ comes from Lemma 3: there exists a constant $C_i$ such that

$$\Gamma_i(x) \leq C_i h^{2d(1-2/s)}.$$

For independent random sequences, $\Gamma_0(x)$ is the only non-zero term.

$$\Gamma_0(x) = C h^d g_2(x) \int K^2(u)\, du \qquad \text{around the point } x. \tag{21}$$

**Proof of Proposition 11.** In order to check assumption (A7), from the Bernstein inequality for independent bounded variables, we get

$$\mathbb{P}(Z(x) > u) \leq 2\exp\left(-\frac{u^2}{2(\operatorname{Var}(Z(x)) + 2Mu\|K\|_\infty/(3\sqrt{nh^d}))}\right). \tag{22}$$

From (21), $\operatorname{Var}(Z(x))$ is bounded by a constant. Because $K$ is a Lipschitz kernel, $x \mapsto Z(x)$ is a Lipschitz function and

$$|Z(x) - Z(y)| \leq 2\sqrt{nh^d} M h^{-(d+1)} \operatorname{Lip} K \|x - y\|_1.$$

Let $M_K$ be the size of the support of $K$ and for $\delta < M_K h$, let $(B_j)_{j=1,\ldots,\nu}$ be a regular partition of diameter $\delta$ over $B$; denoting by $x_j$ the center of $B_j$, we get

$$\sup_{x \in B_j} |K_i(x) - K_i(x_j)| \leq \frac{c\delta}{h} \mathbb{1}\{|x_j - X_i| \leq 2hR\}.$$

Write $\widetilde{g}(x) - \mathbb{E}\widetilde{g}(x) = \widetilde{g}(x) - \widetilde{g}(x_j) + \mathbb{E}(\widetilde{g}(x_j) - \widetilde{g}(x)) + \widetilde{g}(x_j) - \mathbb{E}\widetilde{g}(x_j)$.
For $x \in B_j$,

$$|\widetilde{g}(x) - \widetilde{g}(x_j)| \leq (c\delta/h)\bar{g}(x_j),$$

where

$$\bar{g}(x) = \frac{1}{nh^d} \sum |\widetilde{Y}_i| \mathbb{1}\{|x - X_i| \leq 2hR)\}$$

so that $|\widetilde{g}(x) - \mathbb{E}\widetilde{g}(x)| \leq (c\delta/h)(\bar{g}(x_j) + \mathbb{E}\bar{g}(x_j)) + |\widetilde{g}(x_j) - \mathbb{E}\widetilde{g}(x_j)|$.

Letting $\bar{Z}(x) = \sqrt{nh^d}|\bar{g}(x) - \mathbb{E}\bar{g}(x)|$ and $Z_j = \sup_{x \in B_j} Z(x)$, we have

$$|Z_j| \leq \frac{c\delta}{h}|\bar{Z}(x_j)| + \frac{2c\delta\sqrt{nh^d}}{h} \mathbb{E}\bar{g}(x_j) + |Z(x_j)|.$$

Note that since $\mathbb{E}\bar{g}(x)$ tends to $\bar{g}_1(x) = \int_{|y|\leq M} |y| f(x,y)\,\mathrm{d}y \leq g_1(x)$ as $h$ tends to 0, $2c\delta \mathbb{E}\bar{g}(x)/h \leq t/(3\sqrt{nh^d})$ holds if $\delta \leq c'th/\sqrt{nh^d}$ for a suitable constant $c' > 0$. This condition holds for the following choice (considered only for large values of $t > 0$):

$$\delta = \frac{Ch}{\sqrt{nh^d}} \quad \text{and} \quad 0 < t_0 \leq t. \tag{23}$$

Let $Z = \sup_{x \in B} Z(x)$. For $t > t_0$,

$$\mathbb{P}(Z > t) = \mathbb{P}\Big(\max_{1 \leq j \leq \nu} Z_j > t\Big)$$

$$\leq \nu \max_{1 \leq j \leq \nu} \mathbb{P}(Z_j > t)$$

$$\leq \nu \max_{1 \leq j \leq \nu} \left\{ \mathbb{P}\Big(|\widetilde{Z}(x_j)| > \frac{t}{3}\Big) + \mathbb{P}\Big(|\bar{Z}(x_j)| > \frac{th}{3c\delta}\Big) \right\}$$

$$\leq \nu \max_{1 \leq j \leq \nu} \left\{ \mathbb{P}\Big(|\widetilde{Z}(x_j)| > \frac{t}{3}\Big) + \mathbb{P}\Big(|\bar{Z}(x_j)| > \frac{t}{3}\Big) \right\} \tag{24}$$

$$\leq 4\nu \exp\Big(-\frac{at^2}{1 + tM(nh^d)^{-1/2}}\Big) \tag{25}$$

for some $a > 0$ (from the relation $h/\delta \to \infty$, we assume that $h/c\delta \geq 1$ in order to derive relation (24)). Now,

$$\mathbb{E}Z^p \leq T^p + p\int_T^\infty \mathbb{P}(Z > t) t^{p-1}\,\mathrm{d}t.$$

Choose $T = \sqrt{A \log n}$, note that the function $u \mapsto \exp(-\frac{at^2}{1+tu})$ is non-increasing and assume that $M(nh^d)^{-1/2} \leq T^{-1}$. That is, with the choice for $h$

$$M \leq A^{-1/2} h^{-\rho}, \tag{26}$$

we derive

$$\mathbb{E}Z^p \leq T^p + 4p\nu \int_T^\infty t^{p-1} \exp\Big(-\frac{at^2}{1 + t/T}\Big)\,\mathrm{d}t$$

$$\leq T^p + 4p\nu \int_T^\infty t^{p-1} \exp\Big(-\frac{atT}{2}\Big)\,\mathrm{d}t$$

and using the incomplete gamma function expansion for $x > 2p$

$$\int_x^\infty u^{p-1} \mathrm{e}^{-u}\,\mathrm{d}u \leq 2x^{p-1}\mathrm{e}^{-x}$$

and setting $u = atT/2$ in the previous inequality, we obtain

$$\mathbb{E}Z^p \leq T^p + 4p\nu\left(\frac{2}{aT}\right)^p \int_{aT^2/2}^{\infty} u^{p-1}e^{-u}\,du \leq T^p + \frac{8p}{a}\nu T^{p-1}n^{-aA/2}.$$

With $\nu \sim \delta^{-d}$ and relation (23), the second term is negligible with respect to the first one if $A$ is chosen large enough. Then $\mathbb{E}Z^p = \mathcal{O}((\log n)^{p/2})$ and $\|Z\|_p = \sqrt{nh^d}\|\sup_{x\in B}|\widetilde{g}(x) - \mathbb{E}\widetilde{g}(x)|\|_p = \mathcal{O}((\log n)^{1/2})$. We check that conditions (20) and (26) on $M$ are compatible. This holds if $\rho > dp/(s-2p)$.  □

**Proof of Proposition 12.** First, we define a truncation level $M$ satisfying (20) and consider the truncated variable $\widetilde{g}(x)$. A strong coupling argument by Berbee (1979) yields a Bernstein-type inequality. With Theorem 4 of Doukhan [13], we recall, analogously to (22), that there exist some $\theta, \lambda, \mu > 0$ and an event $A_n$ (which does not depend on either $x$ or $u$) with

$$\mathbb{P}((Z(x) > u) \cap A_n) \leq 4\exp\left(-\frac{\lambda u^2}{2(\mathrm{Var}(Z(x)) + 2Mqu\|K\|_\infty/(3\sqrt{nh^d}))}\right), \qquad (27)$$

with $\mathbb{P}(A_n^c) \leq \mu\beta_{q\theta}$. We first check that $\mathrm{Var}(Z(x))$ is bounded by a constant independent of $x$.

$$|\mathrm{Var}(Z(x)) - h^{-d}\Gamma_0(x)| \leq \frac{2}{h^d}\sum_{i=1}^{n-1}\Gamma_i(x). \qquad (28)$$

As $\beta$-mixing conditions do not improve the bound of the variance, we use a strong mixing condition and the relation $\alpha_i \leq \beta_i$, and refer to the section dedicated to strong mixing. Now, inequality (27) is exactly (22) with $qM$ substituted for $M$. Following along the lines of the proof for the independent case, we get

$$\mathbb{P}((Z > t) \cap A_n) \leq 4\nu\exp\left(-\frac{at^2}{1 + tqM(nh^d)^{-1/2}}\right).$$

The condition (20) on $M$ can now be expressed as $qM \leq h^{-\rho}A^{-1/2}$ and the end of the proof remains unchanged. We get

$$\sqrt{nh^d}\Big\|\sup_{x\in B}|\widehat{g}(x) - \mathbb{E}\widehat{g}(x)|\Big\|_p \leq C\sqrt{\frac{\log n}{nh^d}} + \mathbb{E}Z^p\mathbf{1}_{A_n^c}.$$

Then, using the trivial bound $\mathbb{E}Z^p\mathbf{1}_{A_n^c} \leq \|Z\|_\infty^p\mathbb{P}(A_n^c)$,

$$\sqrt{nh^d}\Big\|\sup_{x\in B}|\widehat{g}(x) - \mathbb{E}\widehat{g}(x)|\Big\|_p \leq C\sqrt{\frac{\log n}{nh^d}} + 2M\|K\|_\infty\frac{\sqrt{n}}{h^{d/2}}\mu\beta_{q\theta}. \qquad (29)$$

Defining $q = n^\gamma$, with $0 < \gamma < 1$, compatibility of the inequalities concerning $M$ requires that $n^{-\gamma}h^{-\rho} > h^{-p(\rho+d)/(s-p)}$. This holds if $\gamma = \frac{\rho(s-2p)-pd}{(2\rho+d)(s-p)}$. Choosing $M =$

$Ch^{-(p/(s-p))(\rho+d)}$, the second term of (29) can be expressed as

$$\frac{Mn}{\sqrt{nh^d}}\beta_{q\theta} = C\frac{\log(n)^{(p/(s-p))(\rho+d)}}{\sqrt{nh^d}}n^{1-\beta\gamma+(p/(s-p))(\rho+d)}.$$

This is negligible with respect to the second term as soon as $\beta > \frac{s\rho+(2s-p)d}{\rho(s-2p)-pd}$. □

**Proof of Proposition 13.** As in the independent case, we choose a truncation level $M = (\frac{n}{\log n})^{(\rho+d)p/((2\rho+d)(s-p))}$ satisfying relation (20) and define $Z(x)$ with respect to the truncated process $\widetilde{g}$. The Fuk–Nagaev inequality leads to

$$\mathbb{P}(Z(x) > u) \leq 4\left(1 + \frac{u^2}{16r\operatorname{Var}(Z(x))}\right)^{-r/2} + \frac{16nM\|K\|_\infty}{u\sqrt{nh^d}}\alpha\left(\frac{u\sqrt{nh^d}}{4M\|K\|_\infty r}\right).$$

We bound $\operatorname{Var}(Z(x))$ above by a constant independent of $x$, using the bound of Lemma 3 and $\Gamma_i(x) = |\operatorname{Cov}(\widetilde{Y}_0 K_0, \widetilde{Y}_i K_i)| \leq 2\int_0^{\alpha_i} Q^2_{\widetilde{Y}K}(t)\,\mathrm{d}t \leq 2\alpha_i M\|K\|_\infty$, and by adapting the proof of Proposition 8. Then, following the proof for the independent case, we get

$$\mathbb{P}(Z > t) \leq c_0\nu\left(1 + \frac{t^2}{c_1 r}\right)^{-r/2} + c_2\nu n\left(\frac{M}{t\sqrt{nh^d}}\right)^{1+\alpha} r^\alpha$$

for suitable constants $c_0$, $c_1$ and $c_2$. We choose $T = \sqrt{A\log n}$ and $r = bT^2$. Then, $\sqrt{\frac{\log n}{nh^d}} = (\frac{n}{\log n})^{-\rho/(2\rho+d)}$. With this choice, $\frac{MT}{\sqrt{nh^d}} = \sqrt{A}(\frac{n}{\log n})^{dp-\rho(s-2p)/((2\rho+d)(s-p))}$.

$$\mathbb{E}Z^p \leq T^p + c_0\nu\int_T^\infty t^{p-1}\left(1 + \frac{t^2}{c_1 bT^2}\right)^{-bT^2/2}\mathrm{d}t$$

$$+ c_2 b^\alpha \nu n\left(\frac{MT}{\sqrt{nh^d}}\right)^{1+\alpha} T^{\alpha-1}\int_T^\infty t^{p-2-\alpha}\,\mathrm{d}t.$$

Changing the variable $u = t/T$, the second term equals

$$c_0\nu T^p\int_1^\infty u^{p-1}\left(1 + \frac{u^2}{c_1 b}\right)^{-bT^2/2}\mathrm{d}u \leq c_0\nu(c_1 b)^{bT^2/2}T^p\int_1^\infty u^{p-1-bT^2}\,\mathrm{d}u$$

$$\leq c_0\frac{\nu T^p(c_1 b)^{bT^2/2}}{bT^2 - p + 1}.$$

Setting $b = 1/(ec_1)$, $(c_1 b)^{bT^2/2} \leq n^{-bA/2}$, it is negligible for a suitable choice of $A$. The third term is less than $\frac{A^{(1+\alpha)/2}}{\alpha-p+1}\nu nT^{p-2}(n/\log n)^{(1+\alpha)(dp-\rho(s-2p))/((2\rho+d)(s-p))}$ with $\nu n = n(n/\log n)^{d(\rho+1)/(2\rho+d)}(\log n)^{d/2}$. Thus, if $dp > \rho(s-2p)$ and $1+\alpha > \frac{(2\rho+2d+d\rho)(s-p)}{dp-\rho(s-2p)}$, then the third term is negligible. The end of the proof follows along the lines of the independent case. □

**Proof of Proposition 14.** Choose $M = (n/\log n)^{(\rho+d)p/((2\rho+d)(s-p))}$ and define $Z(x)$ as in the preceding paragraphs. Proposition 8 and Theorem 1 in [15] imply that

$$\mathbb{P}(Z(x) > u) \leq c_0 \exp\left(-\frac{c_1 u^2}{\mathrm{Var}(Z(x)) + (Mu\|K\|_\infty/\sqrt{nh^d})^{b/(b+2)}}\right)$$

for suitable constants $c_0$ and $c_1$. We first check that $\mathrm{Var}(Z(x))$ is bounded by a constant independent of $x$. The generic term in (28) is bounded by using weak dependence and the fact that the function $u \mapsto K((x-u)/h)$ is $C/h$-Lipschitz:

$$\Gamma_i(x) \leq C\left(\frac{M}{h} + \frac{1}{h^2}\right)\lambda(i) \leq C(h^{-1-(\rho+d)p/(s-p)} + h^{-2})\lambda(i).$$

Considering some $0 < \alpha < 1$, up to a constant, the right-hand side of (28) is bounded above by

$$\sum_i h^d \wedge (h^{-(d+b)}\lambda(i)) \leq h^{d(1-2\alpha)-b\alpha}\sum_i \lambda^\alpha(i),$$

where, in the previous relation, $b = 2$ or $b = 1 + (\rho+d)p/(s-p)$, respectively, if $s \geq (\rho+d+1)p$ or $s < (\rho+d+1)p$. Taking $\alpha < d/(2d+b)$ and noting that $\sum_i \lambda^\alpha(i) < \infty$, the corresponding sum is negligible; thus, $\mathrm{Var}\, Z(x) \sim g_2(x)\int K^2(u)\,\mathrm{d}u$.

Following along the same lines as in the independent case, we then get

$$\mathbb{P}(Z > t) \leq 2c_0\nu \exp\left(-\frac{at^2}{1 + (tM(nh^d)^{-1/2})^{b/(b+2)}}\right)$$

for some $a > 0$. Choosing $T = \sqrt{A\log n}$, $MT/\sqrt{nh^d} < 1$, we have

$$\begin{aligned}\mathbb{E}Z^p &\leq T^p + 4p\nu \int_T^\infty t^{p-1}\exp\left(-\frac{at^2}{1 + (t/T)^{b/(b+2)}}\right)\mathrm{d}t \\ &\leq T^p + 4p\nu \int_T^\infty t^{p-1}\exp\left(-\frac{aT^{b/(b+2)}}{2}t^{(b+4)/(b+2)}\right)\mathrm{d}t\end{aligned} \quad (30)$$

and then setting $u = \frac{aT^{b/(b+2)}}{2}t^{(b+4)/(b+2)}$, the second term of (30) is less than

$$4p\nu\left(\frac{2}{aT^{b/(b+2)}}\right)^p \int_{a/2T^2}^\infty u^{(b+2)p/(b+4)-1}\mathrm{e}^{-u}\,\mathrm{d}u$$

$$\leq 4p\nu\left(\frac{2}{a}\right)^{2p/(b+4)-1} T^{(b^2+4b+8)p/((b+4)(b+2))-1} n^{-aA/2}.$$

With $\nu \sim \delta^{-d}$ and relation (23), the second term of (30) is negligible with respect to $T^p$ if $A$ is chosen large enough. The remainder of the proof is the same as in the independent case. $\square$

# Acknowledgements

We thank Eric Moulines and Alexander Tsybakov for their questions and comments which greatly improved this work. We are also grateful to an anonymous referee for suggesting Corollary 1. This work is dedicated to the memory of Gérard Collomb.

1286 *P. Doukhan and G. Lang*